\documentclass{article}
\usepackage{graphicx} % Required for inserting images
\usepackage{textcomp}
\usepackage{listings}
\usepackage[space=true]{accsupp}
\newcommand{\copyablespace}{\BeginAccSupp{method=hex,unicode,ActualText=00A0}\EndAccSupp{}}
\usepackage[T1]{fontenc}
\usepackage[english]{babel}
\usepackage{hyperref}
\usepackage{graphicx}
\usepackage{amsmath}
\usepackage{amssymb}
\usepackage{amsthm}
\usepackage{mathtools}
\usepackage[table]{xcolor}
\usepackage{fancyhdr}
\usepackage{longtable}
\usepackage{multicol}
\usepackage{enumerate}
\usepackage{enumitem}
\setlist[itemize]{leftmargin=5.5mm}
\usepackage{multirow}
\usepackage{color}
\definecolor{darkgreen}{rgb}{0.1, 0.5, 0.2}
\usepackage{caption}
\usepackage{subcaption}
\usepackage{xurl}
\usepackage[utf8]{inputenc}
\usepackage{tikz}
\usepackage{pgfplots}
\pgfplotsset{compat=newest}

\usepackage{framed}
\usepackage{booktabs}
\usepackage{caption}
\usepackage{float}
\usepackage{titlesec}
\usepackage{capt-of}
\usepackage{array}
\usepackage{arydshln}
\usepackage{appendix}
\newtheorem{theorem}{Theorem}
\newtheorem{definition}{Definition}

\newtheorem{proposition}{Proposition}

\newtheorem{lemma}{Lemma}

\title{Evaluability of paired comparison data in stochastic paired comparison models: Necessary and sufficient condition }
\author{L\'aszl\'o Gyarmati,\ Csaba Mih\'alyk\'o, \'Eva Orb\'an-Mih\'alyk\'o, Andr\'as Mih\'alyk\'o}
\date{}
\begin{document}
\date{\today}
\maketitle

\begin{center}
Department of Mathematics, University of Pannonia, 8200 Veszprém, Hungary\\
email: gyarmati.laszlo@phd.mik.uni-pannon.hu\\
mihalyko.csaba@mik.uni-pannon.hu\\
orban.eva@mik.uni-pannon.hu\\
mihalykoandras@gmail.com\\
\end{center}
%Ezt azért erősen át kellene írni még
\begin{abstract}
In this paper, paired comparison models with stochastic background are investigated. We focus on the models that allow three options for choice. We  estimate all parameters, the strength of the objects and the boundaries of equal decision, by maximum likelihood method. The existence and uniqueness of the estimator are key issues of the evaluation. Although a necessary and sufficient condition for the general case of three options has not been known until now, there are some different sufficient conditions that are formulated in the literature. In this paper, we provide a necessary and sufficient condition for the existence of a maximum and the uniqueness of the argument that maximizes the value, i.e. for the evaluability of the data in models of these types. By computer simulation, we present the efficiency of the condition, comparing it to the previously known sufficient conditions.\noindent
\end{abstract}
\section{Introduction} \label{sec:intro}
Comparisons in pairs are often used in various fields, for example in decision making \cite{MCDM,GROUP}, in analysis of effects \cite{EFFECTS}, in marketing \cite{MARKETING}, in sports \cite{BAKER,SZAD}. There is a very impressive list of applications in publication \cite{VAIDYA}, with more than 150 references.  Paired comparisons are particularly useful when the objects to be evaluated are difficult to relate to the values of a scale but can be compared to each other to determine which of the two is preferred. As the results of the decisions are relations instead of numbers, these types of data require special evaluation methods.\\

One of the main branches of the methods is based on pairwise comparison matrices (PCM), and it is connected to the name of Saaty \cite{SAATY,SAATYII}. Its original method, the Analytic Hierarchy Process (AHP) was generalized in many aspects: the construction of the matrices \cite{KOU}, the evaluation methods \cite{LLSM,Dijkstra}, the elimination of the requirement of complete comparisons \cite{BOZ, CSATO}, and so on. A great amount of publications deal with measuring the inconsistency of the preferences \cite{KAZ,BRUNELLI,AGOSTON,BRUNELLIII}, reduction of inconsistency \cite{ABEL, AGUARON}, determination of the optimal comparison structures \cite{SZAD2, CSATOII}, aggregation methods \cite{AMENTA,DULEBA,SZAD3} and with further issues.\\

Another branch of the paired comparison models has stochastic background. The original idea, the primordial thought, first appears in Thurstone's publications \cite{THURSTONE,THURSTONEII}; therefore, we use the term 'Thurstone-motivated models' when referring to these models. Thurstone envisioned latent probabilistic variables behind the objects to be compared and assumed that decisions are made based on their actual differences. First, two options were allowed for the decisions. Thurstone assumed Gauss distribution, later Bradley and Terry dealt with logistic distribution for the differences \cite{BRADLEY_TERRY}. Another distribution can also be allowed \cite{STERN}, and in \cite{MOEII} it was proved that a large set of distributions have the same properties from axiomatic prospects. The number of options was increased from 2 to 3 \cite{GLENN,RAO}, allowing equal decisions as well. The case of more than three options is considered in \cite{AGRESTI} applying least squares parameter estimation methods and in \cite{MOEI,MOEII} applying maximum likelihood estimation. \\
Maximum likelihood estimation is a very popular, widely applicable parameter estimation method based on optimization. Since the function to be optimized is very complicated, the optimization must be performed numerically. Nevertheless, it is an essential problem whether the maximum value is attained and its argument is unique. Otherwise, the estimated values do not exist and we cannot use them to approximate the strengths or determine the ranking. As the maximum likelihood method is used in the case of Bradley-Terry models, there are numerical methods elaborated and analyzed for this case \cite{HUNTER}. In the case of the two-option Bradley-Terry model, Ford in \cite{FORD} formulated the condition for the existence and uniqueness of the maximizer. Using the special form of the likelihood function, he was also able to prove the existence of the maximum value and uniqueness of the argument. For the three-option Bradley-Terry model, we could not find appropriate results. In a modified three-option model, in the Davidson's model \cite{DAVIDSON}, the author formulated the same condition supplemented by at least one equal decision, and outlined the proof. This paper is the basis of the supplements by $\epsilon$ perturbations made by Conner \cite{CONNER} and later by Yan \cite{YAN}. Yan has made the proof even for further parameters, for example for characterizing the home-field advantage as well. In the case of the three-option Bradley-Terry model, he formulated an equivalence theorem concerning the modified data set. Investigating the problem of the existence and uniqueness of MLE it turned out that Conner and Yan have made unnecessarily many additions: the paper \cite{AXIOMS} contains a much general set of sufficient condition, than the Davidson's condition. But this is not the best condition that can be set up. The authors present an example when the condition in \cite{AXIOMS} is not satisfied, but the maximizer exists and unique. This gap is fulfilled in this paper: using theorems from analysis and graph theory, we prove what are the conditions that are necessary and sufficient for existence and uniqueness for a wide set of three-option stochastic models.\\
In \cite{HAN}, the authors investigate a wide class of extensions of Bradley-Terry models, mainly from the perspective of the asymptotic properties of the maximum likelihood estimation. In the case of three-option models, when considering the problem of existence and uniqueness, they assume the parameter that defines the endpoints of the interval belonging to the equal decision as a predefined constant. This greatly simplifies the proof, but in some cases it cannot be assumed. We do not suppose this; the interval endpoint parameter is also estimated from the data set. \\
The structure of this paper is as follows: in Section \ref{sec:models}
we go through the models, emphasizing their similarities and differences. In Section \ref{sec:cond} we summarize the conditions of evaluability: in \ref{subsec:prev} the previously known conditions, in \ref{subsec:new} the new result, i.e. a necessary and sufficient condition, which characterize the data set from the aspect of evaluability. The proof can be found in Section \ref{sec:app}, Appendix. The proof of sufficiency is a development of the previously used lines of reasoning. In the proof of necessity, graph-theoretical and analytical considerations are connected in an innovative way. In Section \ref{sec:comp} we demonstrate the efficiency of the new set of conditions by computer simulation: we compare the fulfillment rate of the previously known sufficient conditions related to the new finding. Finally, in Section \ref{sec:summary}, we end the paper with a short summary. 

\section{The investigated models}\label{sec:models} \subsection{Thurstone motivated models (THMMs)}\label{subsec:TH}
Let us denote the objects to be evaluated by the numbers $1,2,...,n$. The~Thurstone motivated models have a stochastic background; it is assumed that the current performances of the objects are random variables denoted by $\xi _{i}$, for $i=1,2,...,n$, with the expectations E($\xi _{i}$)=$m_i$. These expectations are the expected strengths of the object $i$. Decisions comparing two objects $i$ and $j$ are related to the difference between these latent random variables, i.e.,~$\xi _{i}$-$\xi _{j}$. We can separate the expectations as follows:
\begin{equation}\label{eq:diff}
\xi _{i}-\xi _{j}=m_i-m_j+\eta_{i,j}
\end{equation}
where $\eta_{i,j}$ are independent identically distributed random variables with the common cumulative distribution function $F$. $F$ is supposed to be three times continuously differentiable, $0<F<1$. The probability density function of F is $0<f$, which is symmetric to zero (i.e.,~$f(-x)=f(x)$, $x \in \mathbb{R}$).  Moreover, the logarithm of $f$ is strictly concave. The set of these c.d.f.-s is denoted by $\mathbb{F}$. Note that logistic distribution and Gauss distribution belong to this set, among many other distributions. 
 If~$\eta_{i,j}$ are Gaussian distributed, then we call it the Thurstone model. If~the distribution $F$ is logistic, i.e.
\begin{equation} \label{eq:logistic}
F(x)=\frac{1}{1+exp(-x)}, x \in \mathbb{R}, 
\end{equation}
we speak about the Bradley--Terry model. These models were originally defined for two-option choices, `worse' and `better', later being generalized for three options: `worse', `better', and `equal'. 
\subsubsection{Two-option models (THMM2)}\label{subsubsec:twoopts}
\begin{figure}[H]
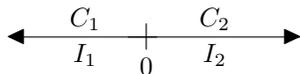

%\hspace{-20mm}\vspace{-20pt}
\setlength{\fboxsep}{0pt}%
\setlength{\fboxrule}{0pt}%
\begin{center}
\[
{\hspace{-23mm}\blacktriangleleft}\hspace*{-0.5038pc}\dfrac{\hspace*{0.75cm}C_{1}%
\hspace*{0.75cm}}{I_{1}}\hspace*{-1.0582pc}%
\begin{array}
[c]{c}%
\\

\mid\\
0%
\end{array}
\hspace*{-1.312pc}\frac{\hspace*{0.95cm}C_{2}\hspace*{0.95cm}}{I_{2}}%
\hspace*{-1.2169pc}%
\hspace*{0.5pc}{\blacktriangleright}%
\
\]
\end{center}
\caption{The options and the intervals belonging to them in a two-option model\label{INTFIG2}.}
\end{figure}

In~the case of the two-option model, the differences are compared to zero: the decision `better'/`worse' indicates whether the difference \eqref{eq:diff} is positive/negative, respectively. The set of real numbers $\mathbb{R}$ is divided into two disjoint parts by the point zero. If the difference of the latent random variables ~$\xi _{i}$-$\xi _{j}$ is negative, then the object $i$ is `worse' than object $j$, while if ~$0 \leq\xi _{i}$-$\xi _{j}$, then $i$ is `better' than $j$. Figure \ref{INTFIG2} represents the intervals of the set of real numbers belonging to the decisions.
If only two options are allowed, then  the probabilities of `worse' and `better' can be expressed as
\begin{equation}\label{eq:veszit2}
p_{i,j,1}=P(\xi_i -\xi_j<0)=F(0-(m_i-m_j))
\end{equation}
and
\begin{equation}\label{eq:gyoz2}
p_{i,j,2}=P(0<=\xi_i -\xi_j)=1-F(0-(m_i-m_j)).
\end{equation}

\vspace{0.5cm}
The data can be included in a three-dimensional ($n \times n \times 2$) data matrix $A$. Its elements $A_{i,j,k}$ represent the number of comparisons in which decision $C_k$ $(k=1,2)$ is the outcome, when object $i$ and $j$ are compared. In detail, $A_{i,j,1}$ is the number of comparisons when object $i$ is `worse' than object $j$, while $A_{i,j,2}$ is the number of comparisons when object $i$ is `better' than object $j$. Of course, $A_{i,j,1}=A_{j,i,2}$.
Assuming independent decisions, the likelihood function is
\begin{equation}\label{eq:LIKL2}
L(A|m_1,...,m_n)=\prod_{k=1}^{2}\prod_{i=1}^{n-1}\prod_{j=i+1}^{n}p_{i,j,k}^{A_{i,j,k}}.
\end{equation}
The log-likelihood function is the logarithm of the above,
\begin{equation}\label{eq:LOGLIKL2}
logL(A|m_1,m_2,...,m_n)=\sum_{k=1}^{2}\sum_{i=1}^{n-1}\sum_{j=i+1}^{n}{A_{i,j,k}\cdot log(p_{i,j,k})}=\sum_{k=1}^{2}\sum_{  i,j=1, i \neq j}^n {0.5 \cdot A_{i,j,k}\cdot log(p_{i,j,k})}.
\end{equation}
One can see that the likelihood function (\ref{eq:LIKL2}), and also the log-likelihood function (\ref{eq:LOGLIKL2}), depend only on the differences of the parameters $m_i$, hence one parameter (or the sum of all parameters) can be fixed. Let us fix $m_1=0$. The maximum likelihood estimation of the parameters $\underline m=(m_1,...,m_n)$ is the argument at the maximum value of (\ref{eq:LIKL2}), or equivalently, of (\ref{eq:LOGLIKL2}); that is
\begin{equation}\label{MLE2}
    \widehat{\underline{m}}=
    \underset{\underline{m}\in \mathbb{R}^{n}, m_1=0}{ \arg \max } \text{ }logL(A|\underline{m}).
\end{equation}
If the distributions of the differences are logistic, 
then the models can be written in the following equivalent form. Introducing 
notations 
\begin{equation} \label{eq:pi}
\pi_ i=\frac{e^{m_i}}{\sum_{l=1}^n e^{m_l}},
\end{equation}
then
\begin{equation}\label{eq:veszit2pi}
p_{i,j,1}=\frac{\pi_j}{\pi_i+\pi_j},
\end{equation}
\begin{equation}\label{eq:gyoz2pi}
p_{i,j,2}=\frac{\pi_i}{\pi_i+\pi_j}
\end{equation}
and the likelihood function can be expressed as
\begin{equation} \label{eq:LIKL2pi}
L(A|\pi_1,...,\pi_n)=\prod_{i=1}^{n-1}\prod_{j=i+1}^{n}(\frac{\pi_j}{\pi_i+\pi_j})^{A_{i,j,1}}(\frac{\pi_j}{\pi_i+\pi_j})^{A_{i,j,2}}
\end{equation}
which (or its logarithm) has to be maximized under the condition 
\begin{equation} \label{eq:feltpi2}
    0<\pi_i, i=1,2,...,n \text{ and } \sum_{i=1}^n\pi_i=1.
\end{equation}
This is the usual form of the two-option Bradley-Terry model  \cite{BRADLEY_TERRY}. We will abbreviate it by BT2 from now on. We note that BT2 is a special case of THMM2.

\subsubsection{Three-option Thurstone motivated models (THMM3)} \label{subsubsec:3opts}
\begin{figure}[H]
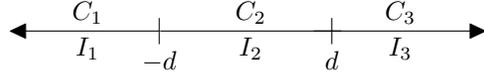

%\hspace{-16mm}\vspace{-16pt}
\setlength{\fboxsep}{0pt}%
\setlength{\fboxrule}{0pt}%
\begin{center}
\[
{\hspace{-23mm}\blacktriangleleft}\hspace*{-0.5038pc}\dfrac{\hspace*{0.75cm}C_{1}%
\hspace*{0.75cm}}{I_{1}}\hspace*{-1.0582pc}%
\begin{array}
[c]{c}%
\\

\mid\\
-d%
\end{array}
\hspace*{-1.0282pc}\dfrac{\hspace*{1cm}C_{2}\hspace*{1cm}}{I_{2}}%
\hspace*{-1.0324pc}%
\begin{array}
[c]{c}%
\\
\mid\\
d%
\end{array}
\hspace*{-1.312pc}\frac{\hspace*{0.95cm}C_{3}\hspace*{0.95cm}}{I_{3}}%
\hspace*{-1.2169pc}%
\hspace*{0.5pc}{\blacktriangleright}%
\
\]
\end{center}
\caption{The options and the intervals belonging to them in a three-option model\label{INTFIG3}.}
\end{figure}

In the case of the two-option model, no parameter is needed to assign the endpoints of the intervals; the number~0 serves for this purpose. By increasing~the number of options to three, however, we need one additional, positive parameter; denoted by $0<d$. The values of the differences are compared to $d$~and/or~$-d$. If the difference is within the interval $[-d,d]$, we can consider the result of the comparison as a new, `equal' decision. Therefore, the set of real numbers $\mathbb{R}$ is divided into three disjoint sub-intervals $I_1$, $I_2$, $I_3$. $I_1=(-\infty,-d)$, $I_2=[-d,d]$,  $I_3=(d,\infty)$. Decisions have a kind of symmetry, i.e.,~if $i$ is `better' than $j$, then $j$ is `worse' than $i$. The figure belonging to the model allowing $3$ options can be seen in Figure~\ref{INTFIG3}.
   
The data can be represented by a three-dimensional ($n \times n \times 3$) data matrix $A$. Its elements $A_{i,j,k}$ represent the numbers of comparisons in which decision $C_k$, $k=1,2,3$ is the outcome when we compare objects $i$ and $j$. The probability that the difference between the random variables belonging to the objects $i$ and $j$ is in the interval $I_k$ can be expressed as follows:
\begin{equation}\label{eq:veszit3}
p_{i,j,1}=P(\xi_i -\xi_j \in I_1)=F(-d-(m_i-m_j)),
\end{equation}
\begin{equation}\label{eq:kozb3}
p_{i,j,2}=P(\xi_i -\xi_j \in I_2)=F(d-(m_i-m_j))-F(-d-(m_i-m_j))
\end{equation}
and
\begin{equation}\label{eq:gyoz3}
p_{i,j,3}=P(\xi_i -\xi_j \in I_3)=1-F(d-(m_i-m_j)).
\end{equation}

The likelihood function is
\begin{equation}\label{eq:LIKL3}
L(A|m_1,...,m_n,d)=\prod_{k=1}^{3}\prod_{i=1}^{n-1}\prod_{j=i+1}^{n}p_{i,j,k}^{A_{i,j,k}}
\end{equation}
and its logarithm 
\begin{equation}\label{eq:LOGLIKL3}
logL(A|m_1,m_2,...,m_n,d)=\sum_{k=1}^{3}\sum_{i=1}^{n-1}\sum_{j=i+1}^{n}A_{i,j,k}\cdot log(p_{i,j,k})=\sum_{k=1}^{3}\sum_{i=1}^{n}\sum_{j=1,j \neq i}^{n} 0.5 \cdot A_{i,j,k} \cdot log(p_{i,j,k}).
\end{equation}
We note that the probabilities \eqref{eq:veszit3}, \eqref{eq:kozb3}  and \eqref{eq:gyoz3} depend only on the differences between the expectations; therefore, one coordinate of the parameter vector $\underline{m}$, for~example $m_1$, can be fixed at zero or the constraint $\sum_{i=1}^{n} m_i=0$ can be assumed. 
The maximum likelihood estimation of the parameters $\underline m=(m_1,...,m_n)$ and $d$ is the argument at the maximal value of (\ref{eq:LOGLIKL3}); that is
\begin{equation}\label{MLE3}
    (\widehat{\underline{m}},\widehat{d})=
    \underset{\underline{m}\in \mathbb{R}^{n},m_1=0,0<d}{ \arg \max } \text{ }logL(A|\underline{m},d).
\end{equation}

Naturally, the~maximal value is not always necessarily attained or the argument is not unique. Nevertheless, some conditions for the data can guarantee the existence of a maximum and uniqueness of its argument. \\
Supposing logistic distributed c.d.f. $F$  (\ref{eq:logistic}), we write the model in the usual form of three-option Bradley-Terry model \cite{RAO}. Introduce notations 
\begin{equation} \label{eq:pi3}
\pi_ i=\frac{e^{m_i}}{\sum_{l=1}^n e^{m_l}},
\end{equation}
and
\begin{equation} \label{eq:d}
\theta=e^d. 
\end{equation}
It is easy to see that 
\begin{equation}\label{eq:veszit3pi}
p_{i,j,1}=\frac{\pi_j}{\theta\cdot\pi_i+\pi_j},
\end{equation}
\begin{equation}\label{eq:kozbpi}
p_{i,j,2}=\frac{(\theta^2-1)\cdot\pi_i\cdot\pi_j}{(\theta\cdot\pi_i+\pi_j)\cdot(\pi_i+\theta\cdot\pi_j)},
\end{equation}

\begin{equation}\label{eq:gyoz3pi}
p_{i,j,3}=\frac{\pi_i}{\pi_i+\theta\cdot\pi_j}.
\end{equation}
The likelihood function is
\begin{equation}\label{eq:LIKL3pi}
L(A|\pi_1,...,\pi_n,\theta)=\prod_{i=1}^{n-1}\prod_{j=i+1}^{n}(\frac{\pi_j}{\theta\cdot\pi_i+\pi_j})^{A_{i,j,1}}\cdot
(\frac{(\theta^2-1)\cdot\pi_i\cdot\pi_j}{(\theta\cdot\pi_i+\pi_j)\cdot(\pi_i+\theta\cdot\pi_j)})^{A_{i,j,2}}\cdot
(\frac{\pi_i}{\pi_i+\theta\cdot\pi_j})^{A_{i,j,3}}
\end{equation}
and its logarithm
\begin{equation}\label{eq:LOGLIKL3pi}
\begin{split}
logL(A|\pi_1,...,\pi_n,\theta)=\sum_{i=1}^{n}\sum_{j=1,j \neq i}^{n}{0.5  \cdot (A_{i,j,1}}\cdot \log(\frac{\pi_j}{\theta\cdot\pi_i+\pi_j})+\\
+{A_{i,j,2}}\cdot\log(\frac{(\theta^2-1)\cdot\pi_i\cdot\pi_j}{(\theta\cdot\pi_i+\pi_j)\cdot(\pi_i+\theta\cdot\pi_j)})
+{A_{i,j,3}}\cdot\log(\frac{\pi_i}{\pi_i+\theta\cdot\pi_j})).
\end{split}
\end{equation}

Recalling (\ref{eq:pi3}) and (\ref{eq:d}), we can see that 
\begin{equation} \label{eq:feltpi3}
0<\pi_i,i=1,...,n,\text{   } 
1<\theta \text{  and   }
\sum_{i=1}^n\pi_i=1
\end{equation}
and the maximum likelihood estimation of the parameters $\underline{\pi}$ and $\theta$ is the argument of maximum value of the likelihood function under the set of conditions (\ref{eq:feltpi3}.) This is the usual form of the three-option Bradley-Terry model. Let us abbreviate it by BT3.
It is easy to see that if $\theta=1$ (that is, $d=0$), we get back to BT2.\\
Finally, we note that the parameters $\pi_i$ can be multiplied by a positive constant without  changing the likelihood/log-likelihood function. Therefore, the optimization can also be performed under the condition $\pi_1=1.$
\subsection{Davidson's model (D3)} \label{sec:DAVIDSON}
In the case of BT2, the Luce's choice axiom is satisfied \cite{LUCE}, but in the case of BT3, it is not. However, to ensure that Luce's axiom is satisfied in the model, Davidson introduced a modified three-option model as follows \cite{DAVIDSON}.\\
Let us consider a vector $\underline{\pi}=(\pi_1,...,\pi_n)$, $0<\pi_i$ satisfying $\sum_{i=1}^n\pi_i=1$ and $0<\nu$. It is assumed that 
\begin{equation} \label{eq:veszitD}
p_{i,j,1}=P(i\text{ is `worse' than }j)=\frac{\pi_j}{\pi_i+\pi_j+\nu \sqrt{\pi_i\*\pi_j}}
\end{equation}
\begin{equation} \label{eq:egyformaD}
p_{i,j,2}=P(i\text{ is `equal' to }j)=\frac{\nu \sqrt{\pi_i\*\pi_j}}{\pi_i+\pi_j+\nu \sqrt{\pi_i\*\pi_j}}
\end{equation}
\begin{equation} \label{eq:gyozD}
p_{i,j,3}=P(i\text{ is `better' than }j)=\frac{\pi_i}{\pi_i+\pi_j+\nu \sqrt{\pi_i\*\pi_j}}
\end{equation}
with some 
\begin{equation} \label{felt2}
0<\nu.
\end{equation}
It is easy to see that by using $\nu=0$ instead, we would get back to BT2. \\
One can check that, due to \eqref{eq:veszitD},\eqref{eq:egyformaD}, and  \eqref{eq:gyozD},
\begin{equation} \label{eq:mertani}
\frac{p_{i,j,2}}{\sqrt{p_{i,j,1}\cdot p_{i,j,3}}}=\nu.
\end{equation}
In this model, the Luce's choice axiom \cite{LUCE} is fulfilled.\\
The likelihood function can be written as follows: 
\begin{equation} \label{eq:likelihood}
L(A|(\pi_1,...,\pi_n,\nu))=\prod_{k=1}^{3}\prod_{i=1}^{n-1} \prod_{j=i+1}^{n} (p_{i,j,k})^{A{i,j,k}}.
\end{equation}
The log-likelihood function is 
\begin{equation} \label{eq:loglikelihood}
log L(A|(\pi_1,...,\pi_n,\nu))=\sum_{k=1}^{3}\sum_{i=1}^{n-1} \sum_{j=i+1}^{n} A_{i,j,k} \cdot log \hspace{0.1cm} (p_{i,j,k})=\sum_{k=1}^{3}\sum_{i=1}^{n} \sum_{j=1 j \neq i}^{n} 0.5 \cdot (A_{i,j,k} \cdot log \hspace{0.1cm} (p_{i,j,k})).
\end{equation}
The maximum likelihood estimation of the parameters ($\widehat{B}$) is the argument of the maximal value of the function \eqref{eq:loglikelihood}, i.e. 
\begin{equation}
\widehat{B}=(\widehat{\underline{\pi }},\widehat{\nu })=\underset{\underline{\pi }>0,\nu
>0,\sum_{i=1}^{n}{\pi_i}=1}{\arg \max } log L(A|(\pi _{1},...,\pi _{n},\nu
))
\end{equation} \label{NEM_ESNEK_EGYBE}
The reader may notice that the parameters $\pi_i$ can be multiplied by a positive constant without changing \eqref{eq:loglikelihood}; therefore, the optimization can also be performed under the condition $\pi_1=1.$\\

After such a detailed presentation of the models, let us prove that D3 does not belong to the set of THMM3, for the sake of completeness.\\

Suppose that there is a cumulative distribution function $F$, $F(-x)=1-F(x), x \in \mathbb{R}$, for which THMM3 coincides with D3. If the strength of the object $i$ equals the strength of the object $j$, i.e. $m_i=m_j$, then, of course, $\pi_i=\pi_j=\pi$. Comparing such objects $i$ and $j$ and, writing the probabilities in both models
\begin{equation}
\begin{split}
p_{i,j,1}=P(\xi_i -\xi_j<-d)=F(-d)=\frac{\pi}{\pi+\pi+\nu\cdot\sqrt{\pi\cdot\pi}}=\\
=1-\frac{\pi}{\pi+\pi+\nu\cdot\sqrt{\pi\cdot\pi}}=\frac{\pi+\nu\cdot\sqrt{\pi\cdot\pi}}{\pi+\pi+\nu\cdot\sqrt{\pi\cdot\pi}}=1-F(d),
\end{split}
\end{equation}
which can be satisfied if and only if $\nu=0$. This contradicts $0<\nu$. Therefore, D3 does not coincide with any Thurstone motivated model.
\section{Conditions of the existence and uniqueness of the maximizer}\label{sec:cond}
In this section we present sufficient and necessary conditions to the various models introduced previously.
\subsection{Previously known sufficient conditions}\label{subsec:prev}
\subsubsection{Two-option model} \label{subsubsec:two-option}
In the case of BT2, Ford formulated the following necessary and sufficient condition for the existence and uniqueness of MLE:
for any partition $S$ and $\overline{S}$, $S\cup\overline{S}={1,2,...,n}$, $S\cap\overline{S}=\emptyset$, there are at least two elements $i$ and $j$, $i\in S$ and $j\in \overline{S}$, for which $0<A_{i,j,2}$.
The statement is proved only for BT2 in \cite{FORD}, using the direct form of the likelihood function substituting (\ref{eq:veszit2pi}) and (\ref{eq:gyoz2pi}). The statement was generalized for a wide class of THMM2, that is, for such cumulative distribution functions which have strictly log-concave density function \cite{AXIOMS}.
The statement can be expressed by the help of graphs as well.
\begin{definition}[Graph of comparisons belonging to the data in two-option models]\label{GR2}
Let the elements 1,2,...,n be represented by nodes of a directed graph, denoted by $\mathbf {GR^{(2)}}$. The nodes $i$ and $j$ are connected with a directed edge (directed from i to j) in $GR^{(2)}$, exactly if $0<A_{i,j,2}$.
\end{definition}
The condition formulated by Ford can be expressed as follows:
in THMMs, allowing two options in choice, the MLE (\ref{MLE2}) exists and is unique if and only if $GR^{(2)}$ is \textit{strongly connected} (that is, there is a directed path from any node $i$ to any other node $j$).

\subsubsection{Three-option Thurstone motivated models}\label{subsbsec:THM3C}
To the best knowledge of the authors, the most general set of conditions for the existence and uniqueness of the MLE (\ref{MLE3}) is contained in \cite{AXIOMS} and looks like:\\

A-SC There is at least one index pair $(i,j)$ for which $0<A_{i,j,2}$.\\

B-SC
For any partition $S$ and $\overline{S}$, $S\cup\overline{S}={1,2,...,n}$, $S\cap\overline{S}=\emptyset$, there are at least two nodes $i$ and $j$, $i \in S$ and $j\in \overline{S}$, for which $0<A_{i,j,2}$, or there are two (not necessarily different) pairs $(i_1,j_1)$ and $(i_2,j_2)$, $i_1, i_2 \in S$, $j_1, j_2 \in \overline{S}$, for which $0<A_{i_1,j_1,3}$ and $0<A_{i_2,j_2,1}$.\\

\begin{definition}[Graph of comparisons in 3-option models] \label{GR^{(3)}} 
Let the nodes of a graph represent the elements 1,2,...,n. The nodes $i$ and $j$ can be connected with two types of edges: one of them is a directed edge belonging to 'better' decisions, the other is a non-directed edge belonging to the 'equal' decision. There is a directed 'better' edge from i to j, if $0<A_{i,j,3}$, and there is non-directed 'equal' edge between $i$ and $j$ if $0<A_{i,j,2}$. This graph will be denoted by $ \mathbf{GR^{(3)}}$.

\end{definition}

C-SC
There is a cycle in the graph $GR^{(3)}$ along the directed edges. (This cycle might contain only two nodes.) \\

The authors proved that the set of conditions A-THMM3, B-THMM3 and C-THMM3 is a sufficient set of conditions, but not a necessary one. While A-THMM3 and B-THMM3 are necessary conditions, an example is provided where the maximal value is attained and its argument is unique, even though C-THMM3 does not hold. Therefore, C-THMM3 can be substituted by a weaker condition, which will be given in Subsection \ref{subsec:new} by C-NS.
\subsubsection{Davidson's model} \label{DC}
In \cite{DAVIDSON}, Davidson formulated the following set of conditions to ensure the existence and uniqueness of MLE \eqref{eq:loglikelihood}. Historically, this is the first condition for the existence and uniqueness of MLE in 3-option models. The method of the proof was a direct method based on the partial derivative of the log-likelihood function.\\

A-DAV There is at least one index pair $ (i,j)$
for which  $0<A_{i,j,2}$.\\

\begin{definition}[Directed graph belonging to the 'better' options in three-option Davidson's models]\label{GR3D}
Let the nodes of a directed graph represent the elements 1,2,...,n. The nodes $i$ and $j$ are connected with a directed edge (from i to j) if $0<A_{i,j,3}$. This graph will be denoted by $\mathbf{GR^{(3D)}}$.  
\end{definition}

B-DAV The graph $G^{(3D)}$, given in Definition (\ref{GR3D}), is strongly connected.\\

We note that condition B-DAV is completely parallel to the condition given by Ford and contained in Subsection \ref{subsubsec:two-option}. Nevertheless, in the case of three-option models, it is far from being a necessary condition. 'Equal' decisions play an important role in evaluability.
We emphasize that this set of conditions is only a sufficient but not a necessary set of conditions. It has been proved in \cite{AXIOMS} that the set of conditions A-SC3, B-SC3 and C-SC3 is a more general set of conditions than that formulated by Davidson.
\subsection{The main result: necessary and sufficient condition of data evaluability in THMM3 and D3} \label{subsec:new}
In this subsection, we formulate the set of necessary and sufficient conditions for evaluability, i.e., the conditions that are necessary and sufficient to ensure the log-likelihood function has a maximal value at a unique argument.\\

Let us assume that the cumulative distribution function $F \in \mathbb{F}$.
Let us define the following graph:
\begin{definition}[graph $GR_{dir}^{(3)}$] \label{GRDIR3}
    Let the objects to be evaluated be assigned the nodes of a graph. Let there be a directed edge from i to j (i $\to$ j) if there exists a decision according to which i is considered 'better' than j, that is, if $A_{i,j,3}>0$. Let there be a special, so called bi-directional edge connecting $i$ and $j$ ($i \longleftrightarrow j$), if there exists a decision, which considers i  'equal' to j, that is, if $A_{i,j,2}>0$. This graph is denoted by $\mathbf {GR_{dir}^{(3)}}$. 
\end{definition}
\begin{definition}\label{directed_cycle}
A directed cycle in $GR_{dir}^{(3)}$ is a sequence of nodes $(i_1,i_2,i_3,....,i_l,i_1)$, where any two consecutive nodes are connected either by a directed 'better' edge ($i_k \to i_{k+1}$) or a bi-directional 'equal' edge ($i_k \longleftrightarrow i_{k+1}$). 
\end{definition}
Note that the direction of the directed 'better' edges matter in a directed cycle. 
\begin{theorem} \label{FOT}
Let $F \in \mathbb{F}$. The maximal values of the log-likelihood functions (\ref{eq:LOGLIKL3}) and (\ref{eq:loglikelihood}) exist and the argument is unique under the conditions $m_1=0$, $0<d$ in THMM3 and $0<\underline{\pi },$ $0<\nu$, $\sum_{i=1}^{n}\pi_i=1$ in D3, respectively, if and only if the data matrix $A$ satisfies the following set of conditions:\\

 A-NS There is at least one index pair $(i,j)$ for which $0<A_{i,j,2}$.\\
 
 B-NS For any partition $S$ and $\overline{S}$, $S\cup\overline{S}={1,2,...,n}$, $S\cap\overline{S}=\emptyset$, there is at least one element $i \in S$ \text{ and } $j\in \overline{S}$, for which $0<A_{i,j,2}$, or there are two (not necessarily different) pairs $(i_1,j_1) \text{ and } (i_2,j_2)$, $i_1, i_2 \in S$, $j_1, j_2 \in \overline{S}$, for which $0<A_{i_1,j_1,3}$ and $0<A_{i_2,j_2,1}$.\\
 
 C-NS  There exists a directed cycle in $GR_{dir}^{(3)}$ (see Definitions~\ref{GRDIR3}) ~and~\ref{directed_cycle}), in which the number of the directed 'better' edges exceeds the number of the bi-directional 'equal' edges.

\end{theorem}
The proof of Theorem \ref{FOT} can be found in Appendix  Section \ref{APPENDIX}.

\section{Efficiency of the necessary and sufficient condition}\label{sec:comp}
In this section, we show that the necessary and sufficient condition for evaluability is met much more frequently than the previously known sufficient conditions in several cases. For this, we applied computer simulations. We generated random comparison data and investigated whether Condition C-NS is satisfied or not, while Condition C-SC is. Of course, this depends on the values of the parameters. If there are a large number of 'better' decisions but only a few 'equal' decisions, we can expect cycles consisting of one-directed edges. If there are few 'better' decisions, the chance of this type of cycle is lower.
We investigated the case $n=8$ in THMM3. The expectations were chosen according to an arithmetic sequence with the difference  $h$, i.e. $\underline{m}=(0,h,2h,...,7h)$. The values of $h$ and $d$ are contained in Table \ref{tab:situations}.

We run $10^9$ simulations for each parameter set. We present similar situations to those in \cite{AXIOMS}, but now we can investigate the efficiency compared to the necessary and sufficient condition.

\begin{table}[H] 
\caption{The choice of parameters} 
\centering
\begin{tabular}{|c|c|c|c|c|}
\hline
Situation& \boldmath{$h$}& \boldmath{$d$}&Rate of Judgments ``Equal''& Rate of ``Better-Worse'' Pairs\\\hline
I.&0.1&0.5&large&large\\\hline
II.&0.1&0.1&small&large\\\hline
III.&0.5&0.5&large&small\\\hline
IV.&0.5&0.1&small&small\\\hline
\end{tabular}
\label{tab:situations} 
\end{table}

\begin{figure*}[!ht]
\setlength{\fboxsep}{0pt}%
\setlength{\fboxrule}{0pt}%
\begin{center}
\includegraphics[width=0.845\textwidth]{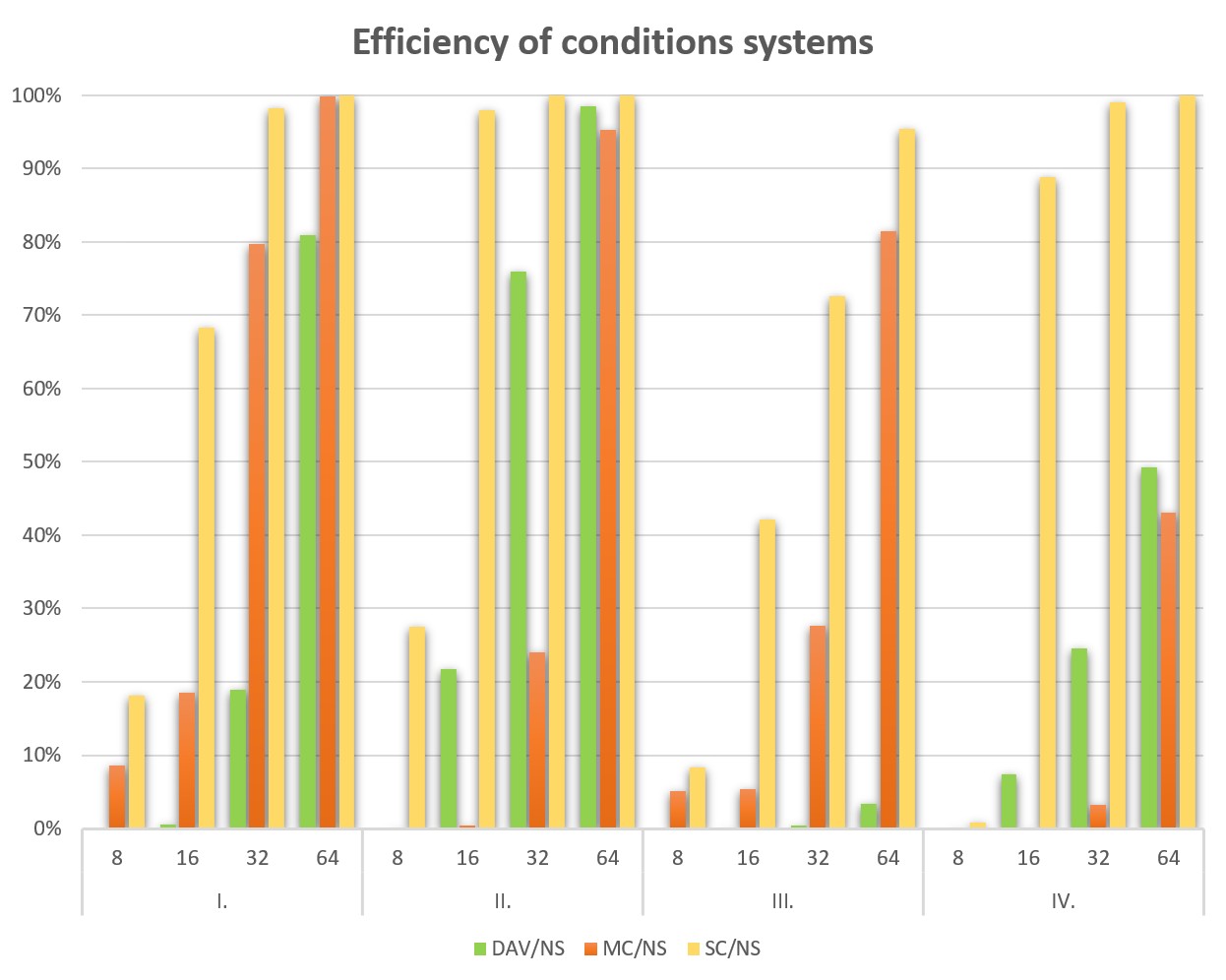}
\end{center}
\caption{Rates of fulfillment of the different condition systems related to the necessary and sufficient conditions.\label{Res}}
\end{figure*}
%Itt nem érthető, hogy melyik modell pontosan mit takar
%ezt alaposan át kell még írni
In Figure \ref{Res} we present the results of the simulations in the case of 8, 16, 32, and 64 comparisons, for every parameter setting of Table \ref{tab:situations}. We randomly selected the edges and also randomly generated the outcomes of decisions (worse/equal/better) according to the probabilities calculated from the parameters. For these random graphs obtained in this way, we checked whether the conditions in \cite{DAVIDSON} (denoted by A-DAV and B-DAV in Subsection \ref{DC}), the conditions in \cite{MOEI} denoted by A-MC, B-MC, C-MC, and the conditions in \cite{AXIOMS} (denoted by A-SC, B-SC and C-SC in Subsection \ref{subsbsec:THM3C}) are satisfied or not. Moreover, we determined the number of such cases when conditions A-NS, B-NS and C-NS in Theorem \ref{FOT} hold.  Taking the ratios of the above quantities by the latest quantity, we get the relative frequencies of the efficiency of the conditions, compared to the necessary and sufficient condition system. The results can be seen in Figure \ref{Res}: when the number of comparisons is small, DAV, MC and SC perform poorly. NS performs much better than DAV, MC and SC. When significantly more comparisons are made (32 or 64), each condition performs better than in the case of few comparisons. If the number of comparison is 8 or 16, there is a big gap between when DAV and MC are fulfilled and when NS is. This is particularly striking in the case of the parameter selection III and IV. SC is more efficient, but in parameter selection III, even with 32 comparisons, we cannot detect more than a quarter of the evaluable data sets. These results show that knowing the necessary and sufficient condition (NS) is very useful, particularly when it is needed to determine whether the data set can be evaluated with a small number of comparisons.
\section{Summary}\label{sec:summary}
In this paper, we studied three-option models with a stochastic background, i.e., three-option Thurstone-motivated models, including the three-option Bradley-Terry model, the generalized Thurstone model, and others, as well as another stochastic model, Davidson's model. We investigated the conditions under which the data can be evaluated by MLE, i.e., the conditions under which the maximum value of the likelihood function exists and is unique. We proved a necessary and sufficient condition system that can characterize the data set from the perspective of evaluability. The proof was developed using widely different areas of mathematics; it is based on the combination of analysis and graph theory.
This condition system proves to be much more effective than the previously known condition systems in certain cases.\\

The authors think that it would be useful to know the necessary and sufficient condition of evaluability in the case of more than three options. Further research is needed to determine this.
\\
 
\section {Appendix}\label{sec:app}

\appendix{}
\section{Appendix\label{APPENDIX}}
\subsection{Sufficiency} \label{subsec:ASUFF}
\subsubsection{THMM3} \label{subsubsec:THMM3II} 
First, we prove that conditions A-NS, B-NS and C-NS in Theorem \ref{FOT}
guarantee the existence of the maximal value of the log-likelihood function  (\ref{eq:LOGLIKL3}) and the uniqueness of its argument. The idea of existence relies on Weierstrass's extreme value theorem: a continuous function on a bounded, closed set attains its minimum and maximum values. For this purpose, we prove that the maximal value of (\ref{eq:LOGLIKL3}) can be searched within a bounded, closed set.\\
Compute $logL(A|\underline{m}=(0,...0),d=1)$=$L_0$. This will serve as a reference point. We note that $L_0<0$. Moreover, each individual term in the sum is negative; therefore, if one of them tends to infinity, the sum itself also tends to infinity.\\
Let $(i,j)$ an index pair for which $0<A_{i,j,2}$ holds. For this $(i,j)$ pair
\begin{equation} \label{d}
A_{i,j,2} \cdot p_{i,j,2}=A_{i,j,2} \cdot (F(d-(m_i-m_j))-F(-d-(m_i-m_j))) \rightarrow - \infty \text{ if } d \rightarrow 0.
\end{equation}
Consequently, condition A-NS guarantees that there exists an 0<$\epsilon$ such that if $d< \epsilon$, then (\ref{eq:LOGLIKL3})<$L_0$. Therefore, the maximal value can be found on the subset
\begin{equation}\label{alsok}
\{d \in \mathbb{R}: \epsilon\leq d\}.
\end{equation}
The upper bound of the parameter $d$ follows from the following reasoning:\\
Take the cycle CY=$(i_1,i_2,...,i_l,i_1)$ in the graph $GR_{dir}^{(3)}$ in Condition C-N). This cycle contains `better' and `equal' edges.
If the edge $i_s,i_{s+1}$ is a 'better' edge, then we can establish that
\begin{equation}
\label{korlgy}
   A_{i_s,i_{s+1},3} \cdot log(1-F(d-(m_{i_s}-m_{i_{s+1}})))\longrightarrow -\infty \text{ supposing  }  d-(m_{i_s}-m_{i_{s+1}})\longrightarrow \infty.
\end{equation} 
Consequently, there exists a value $K_{i_s,i_{s+1}}$ with the following property:
\begin{equation}
\text{if } K_{i_s,i_{s+1}} <d-(m_{i_s}-m_{i_{s+1}}) \text{ then } A_{i_s,i_{s+1},3} \cdot log(1-F(d-(m_{i_s}-m_{i_{s+1}}))) < logL_0.
\end{equation}
Therefore, the maximum has to be within the subset 
\begin{equation}\label{eq:felso_jobb}
d-(m_{i_s}-m_{i_{s+1}}) \leq K_{i_s,i_{s+1}}.
\end{equation}
It means that the maximum can be reached only in such regions where $d-(m_{i_s}-m_{i_{s+1}})$ has an upper bound.\\
If the edge $(i_s,i_{s+1})$ is an `equal' edge, then we can establish that
\begin{equation}
\label{koregyf}
   A_{i_s,i_{s+1},2} \cdot (logF(d-(m_{i_s}-m_{i_{s+1}}))-logF(-d-(m_{i_s}-m_{i_{s+1}})))\longrightarrow -\infty \text{ if }  -d-(m_{i_s}-m_{i_{s+1}})\longrightarrow \infty.
   \end{equation} 
   This can be justified by noticing the fact that 
   \begin{equation} \label{min_plusz}
   d-(m_{i_s}-m_{i_{s+1}})\longrightarrow \infty, \text {if}  -d-(m_{i_s}-m_{i_{s+1}})\longrightarrow \infty.
   \end{equation}
Therefore, there exists a constant value  $K_{i_s,i_{s+1}}$ such that 
\begin{equation}
\begin{aligned}
\text{if } K_{i_s,i_{s+1}}<-d-(m_{i_s}-m_{i_{s+1}}) \text{   then} \hspace{3cm}\\
A_{i_s,i_{s+1},2} \cdot (log(F(d-(m_{i_s}-m_{i_{s+1}}))-F(-d-(m_{i_s}-m_{i_{s+1}}))) < logL_0.
\end{aligned}
\end{equation}
Consequently, the maximum has to be within the subset 
\begin{equation} \label{felso_egyforma}
-d-(m_{i_s}-m_{i_{s+1}}) \leq K_{i_s,i_{s+1}}.
\end{equation}

Sum these inequalities (\ref{eq:felso_jobb}), (\ref{felso_egyforma}) for the edges $(i_1,i_2)$,$(i_2,i_3)$,...$(i_l,i_1)$ belonging to the cycle CY. All expectations $m_{i_s}$ will disappear and we get
\begin{equation} \label{osszead}
d \cdot NU{^{(B)}}+(-d) \cdot NU^{(E)}=d \cdot (NU{^{(B)}}-NU{^{(E)}}) \leq \sum_{(i_s,i_{s+1})\in CY} K_{i_s,i_{s+1}},
\end{equation}
where  $NU^{(B)}$ is the number of the `better'  edges in the cycle CY and $NU^{(E)}$ is the number of the `equal' edges in the cycle CY. By the assumption C-NS, the number of `equal' edges is less than the number of `better' edges, so $0<NU^{(BD)}- NU^{(ED)}$, we can divide by it. We get that we have to search the maximal value of the log-likelihood function within the subset 
\begin{equation} \label{korld}
\{ d \in \mathbb{R}: d\leq \frac{\sum_{(i_s,i_{s+1})\in CY} K_{i_s,i_{s+1}}}{NU^{(B)}- NU^{(E)}} \},
\end{equation}
which is an upper bounded set.
Together with (\ref{alsok}) we could prove that the  region of parameter $d$ can be restricted to a bounded closed set.
From here, the reasoning of the boundedness of the expectations $m_i$ $i=1,2,...,n$ follows the reasoning ST3, ST4 and ST5 in the Appendix of the paper \cite{AXIOMS}, so we do not detail it. As the log-likelihood function (\ref{eq:LOGLIKL3}) is continuous, therefore the maximal value is attained. The uniqueness of the argument is the consequence of the strictly concave property of (\ref{eq:LOGLIKL3}). This property is a part of the theory of log-concave measures \cite{PREKOPA, PREKOPA2}, and is proved in this special case in  the Appendix of the paper \cite{MOE_CJOR2}.
\subsubsection{Davidson's model} \label{subsubsec:Dav_suf}
We strive to make the parametrization of the Davidson model's likelihood function resemble the THMM3 as closely as possible.  We introduce the following form of the parameters: let  
\begin{equation}
\pi_1=1, \pi_i=e^{m_i}, i=1,2,...,n, \text{ and } \nu=e^{\frac{1}{2} \cdot b}.
\end{equation}
We note that condition $\pi_1=1$ corresponds to the condition $m_1=0$.\\
Now, (\ref{eq:veszitD}), (\ref{eq:egyformaD}) and (\ref{eq:gyozD}) can be written in the following forms:
\begin{equation} \label{veszitD2}
    p_{i,j,1}=\frac{1}{e^{(m_i-m_j)}+ e^{\frac{1}{2}(b+m_i-m_j)}+1},
\end{equation}
\begin{equation} \label{egyformaD2}
p_{i,j,2}=\frac{1}{e^{\frac{1}{2}(-b+m_i-m_j)}+e^{\frac{1}{2}(-b+m_j-m_i)}+1},
\end{equation}
and 
\begin{equation} \label{gyozD2}
p_{i,j,3}=\frac{1}{e^{(m_j-m_i)}+e^{\frac{1}{2}(b+m_j-m_i)}+1}.
\end{equation}
Now, $m_1=0$, $m_i \in \mathbb{R}$, $i=2,3,...,n$, $b \in \mathbb{R}$, 
and the log-likelihood function is 
\begin{equation} \label{log-liklD2}
logL(A|m_1,m_2,...,m_n,b)=0.5\sum_{k=1}^{3}\sum_{i\neq j}A_{i,j,k}\cdot log(p_{i,j,k})
\end{equation}
We prove that the maximization can be restricted to a closed bounded set of every parameter $m_i$ and $b$. First we deal with parameter $b$.\\
Let us denote $L_0^{(D)}=logL(A|(0,0,...,0,0)$ in (\ref{log-liklD2}). We use again that if any term in the sum tends to -$\infty$, then the sum itself tends to $-\infty$; therefore, the maximum value cannot be reached in such regions.\\
From (\ref{egyformaD2}), if for an index pair $(i,j)$ $0<A_{i,j,2}$ is satisfied, then 
\begin{equation}
A_{i,j,2}\cdot log (p_{i,j,2}) \longrightarrow -\infty \text{ supposing } b \longrightarrow - \infty,
\end{equation}
as $m_i-m_j\longrightarrow \infty $ and $m_j-m_i\longrightarrow \infty$ cannot hold at the same time. It implies that there exists a value $K_1^{(b)}$ such that if $b<K_1^{(b)}$, then (\ref{log-liklD2})<$L_0^{(D)}$. Therefore, the maximum value of the log-likelihood function can be attained if 
\begin{equation} \label{balso}
K_1^{(b)}\leq b. 
\end{equation}
The upper bound follows from Condition C-NS. It can be proven in a very similar way as in sub-subsection  \ref{subsubsec:THMM3II}.\\
Let us suppose that there exists a cycle $CY^{(D)}$=$(i_1,i_2,...,i_l) \subset GR_{dir}{(3)},$ $(i_s,i_{s+1})$ $s=1,2,...,l$ are edges in the graph $GR_{dir}^{(3)}$ and by assumption C-NS the number of `better' edges is larger than the number of  'equal' edges in the cycle $CY^{(D)}$.\\
If the edge ($i_s,i_{s+1})$ is a `better' edge, then investigate (\ref{gyozD2}).
\begin{equation} \label{felsokd2}
A_{i_s,i_{s+1},3} \cdot log (p_{i_s,i_{s+1},3}) \longrightarrow -\infty \text{ if }  b+m_{i_{s+1}}-m_{i_s} \longrightarrow \infty.
\end{equation}
Therefore, there exists a constant value $K_{i_s,i_{s+1}}^{(U)}$ such that the maximum value of the log-likelihood function is in the region 
\begin{equation} \label{felsok}
b+m_{i_{s+1}}-m_{i_s} \leq K_{i_s,i_{s+1}}^{(U)}
\end{equation}
\\
If the ($i_s,i_{s+1})$ is an 'equal' edge, then we investigate (\ref{egyformaD2}). One can see that
\begin{equation} \label{felsokjd2}
A_{i_s,i_{s+1},2} \cdot log (p_{i_s,i_{s+1},2}) \longrightarrow -\infty \text{ if }  -b +m_{i_{s+1}}-m_{i_s} \longrightarrow \infty.
\end{equation}
This implies that there exists a constant value $K_{i_s,i_{s+1}}^{(U)}$,  the maximal value is in the region
\begin{equation} \label{felsoke2}
-b+m_{i_{s+1}}-m_{i_s} \leq K_{i_s,i_{s+1}}^{(U)}.
\end{equation}
Walking along the cycle $(i_1,i_2,...i_l,i_1)$, we can sum the above inequalities. The parameters $m_s$ disappear and we get  
\begin{equation} \label{korlb}
b \leq \frac{\sum_{(i_s,i_{s+1})\in CY} K_{i_s,i_{s+1}}^{(U)}}{NU^{(BD)}- NU^{(ED)}}.
\end{equation}
where $NU^{(BD)}$ and $NU^{(ED)}$ are the number of `better' edges and 'equal' edges in the cycle $CY^{(D)}$, respectively. 
Inequalities (\ref{balso}) and (\ref{korlb}) guarantee that the maximization with respect to the parameter $b$ can be performed on a closed and bounded subset of $\mathbb{R}$.\\
Let us turn to the parameters $m_i$, $i=1,2,...,n$. We prove that the maximum value of (\ref{log-liklD2}) can be achieved on a closed bounded subset of $\mathbb{R}$ with respect to every index $i
,i=1,2,...,n$. Starting out of $m_1=0$, we prove that every parameter $m_i$ can be restricted to a closed bounded subset.
We show that the lower/upper boundedness is inherited step by step. More precisely, if $K_{i}^{(D)} \leq m_i$  and $0<A_{i,j,2}$, then there exists a value $K_{j}^{(D)}$ such that $L<L_0^{(D)}$ if $m_j< K_{j}^{(D)}$. It means that the maximal value can be achieved on the region $K_{j}^{(D)} \leq m_j$.  As
\begin{equation} \label{expD}
   -b+m_i-m_j \longrightarrow \infty, \text { if } m_j\longrightarrow -\infty,
\end{equation}
consequently, $p_{i,j,2}$ in (\ref{egyformaD2}) tends to 0. Therefore,
\begin{equation}
A_{i,j,2} \cdot log(p_{i,j,2}) \longrightarrow -\infty \text{ supposing }  m_j \longrightarrow -\infty.
\end{equation}
On the other hand, if $0<A_{i,j,1}$, then
\begin{equation}
A_{i,j,1} \cdot log(p_{i,j,1}) \longrightarrow -\infty, 
\end{equation}
as it can be seen from (\ref{veszitD2}), taking into account that
\begin{equation}
    e^{m_i-m_j} \longrightarrow \infty \text{ supposing } m_j \longrightarrow -\infty.
\end{equation} Consequently, the lower bound of the parameter $m_i$ implies a lower bound of the parameter $m_j$.\\
Now, recalling Condition B-NS, starting from $i_1$, it is connected by an `equal' or a `better' edge with an element denoted by $i_2$, this latter edge is directed from $i_2$ to $i_1$. Consequently, $m_{i_2}$ has a lower bound. Defining $S=\{i_1,i_2\}$  and $\overline{S}=\{1,...,n\} \setminus S$
, there exists an `equal' edge or `better' edge (directed from an element $i_3  \text{ in } \overline S$ to an element of $S$ , which guarantees the lower bound for $m_{i_3}$.
This process can be continued until $S$ becomes empty set. \\
Exactly the same can be stated walking along the `better' and 'equal' edges. An upper bound of a parameter $m_i$ implies an upper bound for those parameters which belong to such nodes which are connected with $i$ by an 'equal' or a `better' edge. Starting from the node $1$, $m_1=0$, therefore each parameter $m_j$ has an upper bound in the sense that under this upper bound the value of the log-likelihood function (\ref{log-liklD2}) is larger than $L_0^{(D)}$. Therefore the maximization of the likelihood function can be restricted to a bounded closed region of $\mathbb{R}^{n+1}$, and the log-likelihood function is continuous; therefore the maximal value is attained.\\

The uniqueness of the argument is again the consequence of the strictly concave property of the log-likelihood function (\ref{log-liklD2}). One can check that the Hesse matrices of the functions (\ref{veszitD2}), (\ref{egyformaD2}) and (\ref{gyozD2}), as a function of $x=m_i-m_j$ and $b$ are negative definite for any value of $x$ and $b$; therefore, the sum is also a strictly concave function of all its variables; therefore, the argument of the maximum value is unique.

\subsection{Necessity} \label{subsec:ANEC}
\subsubsection{THMM3} \label{subsubsec:THMM3III}
First, let us investigate condition A-NS in Theorem \ref{FOT}. Recalling formulas (\ref{eq:veszit3}) and (\ref{eq:gyoz3}), the probabilities $p_{i,j,1}$ and $p_{i,j,3}$ are strictly increasing if $d \to 0+$. Due to the lack of $p_{i,j,2}$ terms, the log-likelihood function (\ref{eq:LOGLIKL3}) is also increasing. Consequently, it does not attain its maximum value. We conclude that Condition A-NS is a necessary condition for the function to have a maximal value.\\
\\
Turn to Condition B-NS. 
If Condition B-NS does not hold, then there  exists a partition $S$ and $\overline{S}$ such that there is no `equal' edge and `better' edge from $S$ to $\overline{S}$ and from $\overline{S}$ to $S$. If there exist edges in $GR^{(3)}$ between these subsets, their direction is the same. Suppose that there is a `better' edge from $S$ to $\overline{S}$. Let $i_1,...,i_k \in S,$ $j_1,...,j_{n-k} \in \overline{S}$.\\
Suppose that there is a unique maximizer at $(m_1^{(0)},...,m_n^{(0)})$ and $d^{(0)}.$ In this case, the log-likelihood function contains only $p_{i,j,3}$ terms and not $p_{i,j,1}$ and $p_{i,j,2}$ terms ($i \in S$
and $j \in \overline{S}$). Increase the values of the expectations belonging to the indices in $S$ by 1 and keep the values of the expectations belonging to the indices in $\overline{S}$, i.e. let $m_i^{(1)}=m_i^{(0)}+1$ if $i \in S$ and $m_j^{(1)}=m_j^{(0)}$ if $j \in \overline{S}$. Take the log-likelihood  at $\underline{m}^{(1)}$. In this case, those terms which belong to such index pairs $(i,j)$ for which $i,j \in S$ or $i,j \in \overline{S}$, do not change. Those terms would definitely increase if $i$ and $j$ are in different subsets. Therefore, we could increase the value of the likelihood and the log-likelihood function, which is a contradiction. We note that if there is no edge between $S$ and $\overline{S}$, making the same changes, the value of log-likelihood does not increase but the maximum value is attained at $ \underline{m}^{(0)}$ and $\underline {m}^{(1)}$ too, which contradicts the uniqueness of the argument.\\

Let us now tackle the necessity of Condition C-NS. Suppose again that there exists a unique maximizer at $(m_1,\dots ,m_n)$ and $d$. We shall modify the values of $(m_1,\dots ,m_n)$ and $d$ while increasing the value of the log-likelihood function \eqref{eq:LOGLIKL3} achieving a contradiction. For this, we need to investigate potential ways to modify some values without decreasing the value of \eqref{eq:LOGLIKL3}. The easiest way to do so is by changing only the value of $d$.
\begin{lemma} \label{Increase_only_d}   
For any $x>0$,

\begin{equation} \label{d_ben_no}
F(d-(m_i-m_j))-F(-d-(m_i-m_j))<F(d+x-(m_i-m_j))-F(-d-x-(m_i-m_j)). 
\end{equation}
\end{lemma}
That is to say, if $d$ increases while $m_i$ and $m_j$ remain unchanged, the probability of `equal' decision strictly increases. \\

We can state similar results even if $m_i$ and $m_j$ change alongside with $d$. We can list the potential ways how we can introduce these changes so that the probability cannot change. The results listed in the following lemma are quite intuitive, still, let us state them formally. 
\begin{lemma} \label{Modified_inequalities}

Let us define a new set of values of parameters as follows:
\begin{equation} \label{MOD}
d^{(1)}=d+x, m_i^{(1)}=m_i+y \text{ and } m_j^{(1)}=m_j+z.
\end{equation}

Let us define $p_{i,j,1}^{(1)}, p_{i,j,2}^{(1)}$ and $p_{i,j,3}^{(1)}$ similarly to equations \eqref{eq:veszit3}, \eqref{eq:kozb3} and \eqref{eq:gyoz3} as follows:
\begin{equation}
p_{i,j,1}^{(1)}=F(-d^{(1)}-(m_i^{(1)}-m_j^{(1)})),
\end{equation}
\begin{equation}
p_{i,j,2}^{(1)}=F(d^{(1)}-(m_i^{(1)}-m_j^{(1)}))-F(-d^{(1)}-(m_i^{(1)}-m_j^{(1)}))
\end{equation}

\begin{equation}
p_{i,j,3}^{(1)}=1-F(d^{(1)}-(m_i^{(1)}-m_j^{(1)})).
\end{equation}

The following inequalities hold:
\begin{itemize}

\item Suppose $i$ and $j$ are connected by an `equal' edge, i.e. $0<A_{i,j,2}$,\\ 

\begin{equation}
\text{If  }\left| y-z \right| \leq x, \text{ then  } A_{i,j,2}\cdot \log (p_{i,j,2} )\leq A_{i,j,2} \cdot \log (p_{i,j,2}^{(1)}).
\end{equation}
\begin{equation}
\text{If  }\left| y-z \right|<x, \text{ then  } A_{i,j,2} \cdot \log (p_{i,j,2}) < A_{i,j,2} \cdot \log (p_{i,j,2}^{(1)}).
\end{equation}

\begin{equation}
\text{If  } \left| y-z \right|>x, \text{ then  } A_{i,j,2} \cdot \log (p_{i,j,2}^{(1)}) < A_{i,j,2} \cdot log (p_{i,j,2}).
\end{equation}

\item Suppose that $i$ and $j$ are connected by a `better' edge, i.e. $0<A_{i,j,3}$. Apply the same parameter modification as in (\ref{MOD}).\\
\begin{equation}
    \text{If  } x \leq y-z, \text{ then } A_{i,j,3} \cdot \log (p_{i,j,3}) \leq A_{i,j,3}\cdot \log (p_{i,j,3}^{(1)}).
    \end{equation}
    \begin{equation}
     \text{If  } x < y-z, \text{ then } A_{i,j,3} \cdot \log (p_{i,j,3}) < A_{i,j,3}\cdot \log (p_{i,j,3}^{(1)}).
     \end{equation}
     \begin{equation}
      \text{If  } y-z < x, \text{ then } A_{i,j,3}\cdot \log (p_{i,j,3}^{(1)}) < A_{i,j,3}\cdot \log (p_{i,j,3}).
    \end{equation}
\item An analogous statement can be stated for the case if $i$ and $j$ are connected by a `worse' edge, which is a better edge from $j$ to $i$, as $0<A_{i,j,1}$=$A_{j,i,3}$.
\end{itemize}
\end{lemma}

The proof of Lemma \ref{Modified_inequalities} is a technical exercise; it is left to the reader. 

Lemmas~\ref{Increase_only_d} and \ref{Modified_inequalities} immediately imply the following: 

\begin{lemma}\label{together-change}
    If there is a unique maximizer at $(m_1,\dots ,m_n)$ and $d$, then there are no $(m_1^{(1)},\dots ,m_n^{(1)})$ and $d^{(1)}$ values, with $d^{(1)}-d = x > 0$ and $m_i^{(1)} - m_i = u_i$ for which

    \begin{itemize}
    \item for every $i$ and $j$ which are connected by an `equal' edge, (i.e. $0<A_{i,j,2}$) \text{  } $\left| u_i-u_j \right| \leq x$  
            \item for every $i$ and $j$ which are connected by a `worse' edge, (i.e. $0<A_{i,j,1}=A_{j,i,3}$) \text{  } $ u_i-u_j \leq -x$

        \item for every $i$ and $j$ which are connected by a `better' edge, (i.e. $0<A_{i,j,3}$) \text{  } $u_i-u_j \geq x$

    \end{itemize}
    
\end{lemma}

In essence Lemma~\ref{together-change} describes, if we can change the values of $(m_1,\dots,m_n)$ and $d$ in a such a way that it respects the conditions of Lemmas~\ref{Increase_only_d} and \ref{Modified_inequalities}, we can increase the value of log-likelihood function. Let us now attempt to modify the values in such way. 

Let us form the following weighted directed multigraph $GW$ from the values of the comparison.
\begin{definition} [Weighted graph $GW$]
    Let the objects to be evaluated be assigned the nodes of $GW$. Let there be a directed edge from $i$ to $j$ ($i \to j$) with weight $-1$, if there is a decision according to $i$ is considered 'better' than $j$, that is, if $A_{i,j,3}>0$. Let there be bi-directed edges from $i$ to $j$ and $j$ to $i$ ($i \longleftrightarrow j$) with weight $+1$, if $i$ is considered 'equal' to $j$ according to a decision, that is, if $A_{i,j,2}>0$. 
\end{definition}
In other words, there are -1 or +1 on the weights of the edges of graph $GR_{dir}^{(3)}$ in Definition \ref{GRDIR3}.
Note that there might be directed edges of weight $-1$ and $+1$ thus between $i$ and $j$ both directions.\\
\begin{definition}[length of a path]
The length of a path is the sum of the weights belonging to the edges of the path. A shortest path is a path of minimum length.

\end{definition}
Let us consider an arbitrary vertex, for example node $1$. Let us find the shortest path in $GW$ from $1$ to every other vertex. 
\begin{lemma} 
  Condition B-NS implies that there is a directed path in $GW$ from $1$ to every other vertex.   
\end{lemma}
\begin{proof}
    Suppose that there is a vertex $j$, that cannot be reached from $1$.  Let $S$ be the set of vertices not reachable from $1$. $1\in S\neq \emptyset$, $j \in \overline S \neq \emptyset$. But there is no 'better' or 'equal' edges from $S$ to $\overline{S}$, which contradicts Condition B-NS.
\end{proof}
Now we use the following well-known statement: we can find the shortest path from $1$ to every other vertex in $GW$ if and only if $GW$ is conservatively weighted, that is, if there is no negative circle in $GW$ \cite{FRANK} (a negative circle is a list of edges ${(i_1,i_2,...,i_l,i_1)}$, the sum of their weights is negative). 

\begin{lemma}
    Negation of Condition C-NS is equivalent to the conservative property of the weighted graph $GW$.
\end{lemma}
\begin{proof}
Negation of Condition C-NS expresses that there is no such cycle in the graph $GW$, which contains more 'better' edges than 'equal' edges. This means that along every cycle, the sum of the weights +1  is greater than or equal to the sum of the weights -1; therefore, there is no negative cycle.
\end{proof}

Thus, there is a shortest path from $1$ to every other vertex. Let the shortest distance from $1$ to $i$ be denoted by $u_i$ ($u_1 = 0$ by the conservative weights of $GW$). Notice the following observations:
\begin{lemma}\label{correct_weights}
\
\begin{itemize}
    \item If there is an 'equal' edge between $i$ and $j$, then there is an edge with weight $+1$ between $i\to j$ and $j\to i$. Thus the shortest path from $1$ to $j$ cannot be longer than $u_i+1$ and the shortest path from $1$ to $i$ cannot be longer than $u_j+1$, hence $|u_i-u_j|\leq 1$.
    \item If there is a 'better' edge between $i$ and $j$, then there is an edge with weight $-1$ between $i\to j$. Thus the shortest path from $1$ to $j$ cannot be longer than $u_i-1$. Hence $u_j\leq u_i -1$.
\end{itemize}
    
\end{lemma}
\begin{proposition}
    If there exists a unique maximizer of the likelihood function, then Condition C-NS is held.

\end{proposition}

\begin{proof}

Suppose there exists a unique maximizer of the likelihood function and Condition C-NS does not hold. Let us define $ d^{(1)}= d+1$ and $(m_1^{(1)},\dots ,m_n^{(1)}) = (m_1+u_1,\dots ,m_n+u_n)$. We can see by Lemma~\ref{correct_weights} that $d^{(1)}$ and $(m_1^{(1)},\dots ,m_n^{(1)})$ contradict Lemma~\ref{together-change}. Hence Condition C-NS is necessary to have a unique maximizer. 
\end{proof}

\subsubsection{Davidson's model} \label{subsubsec:D3}
First, let us investigate condition A-NS in Theorem \ref{FOT}. Recalling formulas (\ref{veszitD2}) and (\ref{gyozD2}), we can see that probabilities $p_{i,j,1}$ and $p_{i,j,3}$ are increasing if $b \to -\infty$, and due to $A_{i,j,2}=0$  $i=1,...,n$ $j=1,...,n$ $i \neq j$, $A_{i,j,2} \cdot p_{i,j,2}$ terms are missing from the sum. Therefore, the log-likelihood function \ref{log-liklD2} does not attain its maximal value. We conclude that Condition A-NS is a necessary condition for the function to have a maximal value.\\
Turning to condition B-NS, the same reasoning as in sub-subsection \ref{subsubsec:THMM3III} applies.\\
Now, let us investigate Condition C-NS.
Starting out of the form \ref{veszitD2}, \ref{egyformaD2} and \ref{gyozD2},  
one can easily check that, similarly to \eqref{d_ben_no}  in Lemma (\ref{Increase_only_d}),
\begin{lemma} \label{Increase_only_d_D3}   
For any $x>0$,
\begin{equation} \label{egyformaD3_II}
p_{i,j,2}=\frac{1}{ e^{\frac{1}{2}(-b+m_i-m_j)}+e^{\frac{1}{2}(-b-(m_i-m_j)}+1}<\frac{1}{ e^{\frac{1}{2}(-(b+x)+m_i-m_j)}+e^{\frac{1}{2}(-(b+x)-(m_i-m_j)}+1}=p_{i,j,2}^{(1)}
\end{equation}

\end{lemma}
Again, if $b$ increases while $m_i$ and $m_j$ remain unchanged, the probability of `equal' decision strictly increases. Consequently, Lemma \ref{Modified_inequalities} with its each inequality remains true.
The reader can check the following: 
%------------------
\begin{lemma} \label{Modified_inequalitiesD}

Let us define a new set of values of parameters as follows:
\begin{equation} \label{MOD_D}
b^{(1)}=b+x,\text{  } m_i^{(1)}=m_i+y \text{ and } m_j^{(1)}=m_j+z.
\end{equation}

Let us define $p_{i,j,1}^{(1)}, p_{i,j,2}^{(1)}$ and $p_{i,j,3}^{(1)}$ similarly to equations \eqref{veszitD2}, \eqref{egyformaD2} and \eqref{gyozD2} as follows:
\begin{equation}
 p_{i,j,1}^{(1)}=\frac{1}{e^{(m_i^{(1)}-m_j^{(1)})}+ e^{(b^{(1)}+m_i^{(1)}-m_j^{(1)})}+1}
 \end{equation}
\begin{equation}
p_{i,j,2}^{(1)}=\frac{1}{e^{\frac{1}{2}(-b^{(1)}+m_i^{(1)}-m_j^{(1)})}+e^{\frac{1}{2}(-b^{(1)}+m_j^{(1)}-m_i^{(1)})}+1},
\end{equation}

\begin{equation}
p_{i,j,3}^{(1)}=\frac{1}{e^{(m_j^{(1)}-m_i^{(1)})}+e^{\frac{1}{2}(b^{(1)}+m_j^{(1)}-m_i^{(1)})}+1}.
\end{equation}

The following inequalities hold:
\begin{itemize}

\item Suppose $i$ and $j$ are connected by an `equal' edge, i.e. $0<A_{i,j,2}$,\\ 

\begin{equation}
\text{If  }\left| y-z \right| \leq x, \text{ then  } A_{i,j,2}\cdot \log (p_{i,j,2} )\leq A_{i,j,2} \cdot \log (p_{i,j,2}^{(1)}).
\end{equation}
\begin{equation}
\text{If  }\left| y-z \right|<x, \text{ then  } A_{i,j,2} \cdot \log (p_{i,j,2}) < A_{i,j,2} \cdot \log (p_{i,j,2}^{(1)}).
\end{equation}
\begin{equation}
\text{If  } \left| y-z \right|>x, \text{ then  } A_{i,j,2} \cdot \log (p_{i,j,2}^{(1)}) < A_{i,j,2} \cdot log (p_{i,j,2}).
\end{equation}

\item Suppose that $i$ and $j$ are connected by a `better' edge, i.e. $0<A_{i,j,3}$. Apply the same parameter modification as in (\ref{MOD_D}).\\
\begin{equation}
    \text{If  } x \leq y-z, \text{ then } A_{i,j,3} \cdot \log (p_{i,j,3}) \leq A_{i,j,2}\cdot \log (p_{i,j,3}^{(1)}).
    \end{equation}
    \begin{equation}
     \text{If  } x < y-z, \text{ then } A_{i,j,3} \cdot \log (p_{i,j,3}) < A_{i,j,3}\cdot \log (p_{i,j,3}^{(1)}).
     \end{equation}
     \begin{equation}
      \text{If  } y-z < x, \text{ then } A_{i,j,3}\cdot \log (p_{i,j,3}^{(1)}) < A_{i,j,3}\cdot \log (p_{i,j,3}).
    \end{equation}
\item An analogous statement can be stated for the case if $i$ and $j$ are connected by a `worse' edge, which is a `better' edge from $j$ to $i$ i.e. $0<A_{i,j,1}$=$A_{j,i,3}$.
\end{itemize}
\end{lemma}

%-------------------
This lemma is entirely analogous to Lemma  \ref{Modified_inequalities}. Therefore, all further reasoning  in sub-subsection \ref{subsubsec:THMM3III} applies, demonstrating the necessity of Condition C-NS in the case of D3 as well.
\end{document}